\documentclass{article}

\usepackage{graphicx}
\usepackage{amsmath}
\usepackage{amsthm}
\usepackage{amsfonts}
\usepackage{amssymb}
\usepackage{url}

\newtheorem{thm}{Theorem}[section]
\newtheorem{defn}[thm]{Definition}
\newtheorem{prop}[thm]{Proposition}

\newtheorem{lem}[thm]{Lemma}

\newcommand{\Z}{\mathbb{Z}}

\newcommand{\Gm}{\Gamma}

\newcommand{\bz}{{\bf 0}}
\newcommand{\bo}{{\bf 1}}
\newcommand{\bzT}{{\bf 0}^T}
\newcommand{\boT}{{\bf 1}^T}

\newcommand{\bit}{\begin{itemize}}
\newcommand{\eit}{\end{itemize}}
\newcommand{\ben}{\begin{enumerate}}
\newcommand{\een}{\end{enumerate}}
\newcommand{\beq}{\begin{equation}}
\newcommand{\eeq}{\end{equation}}
\newcommand{\bea}{\begin{eqnarray*}}
\newcommand{\eea}{\end{eqnarray*}}
\newcommand{\bpf}{\begin{proof}}
\newcommand{\epf}{\end{proof}\ms}
\newcommand{\ms}{\medskip}
\newcommand{\noi}{\noindent}

\title{Construction of Directed Strongly Regular Graphs Using Block Matrices}
\author{Francis Adams\thanks{Department of Mathematics, Carleton College, Northfield,
MN 55057, USA (adamsf@carleton.edu).  Research supported by DMS
0750986 (REU-site)} \and Alexandra Gendreau\thanks{Department of
Mathematics, Wellesley College, Wellesley, MA 02481, USA
(agendrea@wellesley.edu). Research supported DMS 0750986 (REU-site)}
\and Oktay Olmez\thanks{Department of Mathematics, Iowa State
University, Ames, IA 50011, USA (oolmez@iastate.edu).} \and Sung
Yell Song\thanks{Department of Mathematics, Iowa State University,
Ames, IA 50011, USA (sysong@iastate.edu).  } }

\begin{document}
\maketitle

\begin{abstract}
The concept of directed strongly regular graphs was introduced by
Duval in ``A Directed Graph Version of Strongly Regular Graphs"
[\textit{Journal of Combinatorial Theory}, Series A 47 (1988) 71 --
100]. Duval also provided several construction methods for directed
strongly regular graphs. We construct several new classes of
directed strongly regular graphs with parameters $\lambda=\mu=t-1$
or $\lambda+1=\mu=t$. The directed strongly regular graphs reported
in this paper are obtained using a block construction of adjacency
matrices of regular tournaments and circulant matrices. We then give
some algebraic and combinatorial interpretation of these graphs in
connection with known directed strongly regular graphs and related
combinatorial structures.
\end{abstract}

\noi{\bf Keywords} Directed Strongly Regular Graphs, Cayley Graphs,
Regular Tournaments, Doubly-Regular Team Tournaments

\ms\noi{\bf AMS Classification}

\ms

\section{Introduction}\label{sintro} 
This paper investigates directed strongly regular graphs and some methods of constructing them from block matrices.  Section 2 introduces the necessary notation and defines directed strongly regular graphs in terms of its adjacency matrix and its combinatorial properties.  Section 3 looks at feasibility condtions of parameter sets established by Duval and some construction methods he provides. We also describe a DSRG construction method not used by Duval, construction using Cayley Graphs of groups.  Some of these Cayley Graph construction motivate the first constructions in Section 4 based on regular tournaments.  After these and other constructions stemming from block matrices, especially regular tournaments, Section 5 investigates when these constructions produce isomorphic graphs, including using different tournaments in one construction and the same tournament in different constructions.  Section 7 is a summarizing list of all our construction methods and the parameter sets satisfied by each, as well as a short list of the first few parameter sets for which we have found new constructions.

\section{Preliminaries}

All graphs considered in this paper will be finite simple graphs; so
our graphs will have no loops or multiple edges. Let $\Gamma$ be a
directed graph with vertex set $V(\Gamma$) and edge set $E(\Gamma$).
For any $x,y \in$ V($\Gamma$), we say that $x$ is adjacent to $y$,
denoted $x\rightarrow y$, if there is an edge from $x$ to $y$. There
may also be an edge from $y$ to $x$, in which case we will say there
is an undirected edge between $x$ and $y$, written as
$x\leftrightarrow y$.  Finally, $x$ is not adjacent to $y$,
signified as $x\nrightarrow y$, if neither $x\rightarrow y$ nor
$x\leftrightarrow y$.

Let $\Gamma$ be a graph with $V(\Gamma)=\{x_1, x_2,\dots, x_n\}$.
The adjacency matrix $A=A(\Gm)$ of $\Gamma$ is an $n\times n$ matrix
whose rows and columns are indexed by the vertices such that
$$A_{ij}=\left \{\begin{array}{ll} 1 & \mbox{if } x_i \mbox{ is adjacent
to } x_j\\ 0 & \mbox{otherwise}\end{array}\right .$$  We will use $\Gamma$ and its adjacencey matrix $A$ interchangeably.

A strongly regular graph with parameters ($n, k, \lambda, \mu$) is
an undirected graph with $n$ vertices whose adjacency matrix $A$
satisfies the following equations: $$A^2=kI + \lambda{A} +
\mu(J-I-A)$$ $$AJ=JA=kJ$$ where $I$ is the identity matrix and $J$
is the all-ones matrix. From the first equation, we see that the
number of paths of length two from a vertex $x$ to another vertex
$y$ is $\lambda$ if $x$ and $y$ are adjacent, $\mu$ if $x$ and $y$
are not adjacent. This second equation means that each vertex has
valency $k$.

The concept of `strong regularity' in the class of directed graphs
is a generalization of that in the class of undirected graphs. Let
$\Gamma$ be a directed graph with its adjacency matrix $A$. The
graph $\Gamma$ is called a \textit{directed strongly regular graph
with parameters} $(n, k, t, \lambda, \mu)$, denoted
DSRG$(n,k,t,\lambda, \mu)$, if $A$ satisfies the following
equations:
$$A^2=tI + \lambda{A} + \mu(J-I-A)$$
$$AJ=JA=kJ$$  Thus, each vertex of DSRG$(n,k,t,\lambda, \mu)$ has
$k$ out-neighbors and $k$ in-neighbors, including $t$ neighbors
counted as both in- and out-neighbors of the vertex. For vertices $x
\neq y$, there are $\lambda$ paths of length two from $x$ to $y$ if
$x \rightarrow y$ and $\mu$ paths of length two if $x \nrightarrow
y$.

A strongly regular graph with parameters $(n,k,\lambda, \mu)$ is
viewed as a DSRG$(n, k, t, \lambda, \mu)$, characterized by a DSRG
with $t=k$. A DSRG with $t=0$ is a graph known as
doubly-regular tournament. An example of a DSRG that is not a
strongly regular graph has the parameters $(8, 3, 2, 1, 1)$ with
adjacency matrix $A$ below. This DSRG is illustrated in Figure 1.
$$A = \left[ \begin{array}{cccccccc}
0 & 1 & 0 & 0 & 1 & 0 & 0 & 1 \\
0 & 0 & 1 & 0 & 1 & 1 & 0 & 0 \\
0 & 0 & 0 & 1 & 0 & 1 & 1 & 0 \\
1 & 0 & 0 & 0 & 0 & 0 & 1 & 1 \\
1 & 1 & 0 & 0 & 0 & 0 & 0 & 1 \\
0 & 1 & 1 & 0 & 1 & 0 & 0 & 0\\
0 & 0 & 1 & 1 & 0 & 1 & 0 & 0\\
1 & 0 & 0 & 1 & 0 & 0 & 1 & 0\\\end{array} \right] $$
\begin{center}
\includegraphics[scale=0.3]{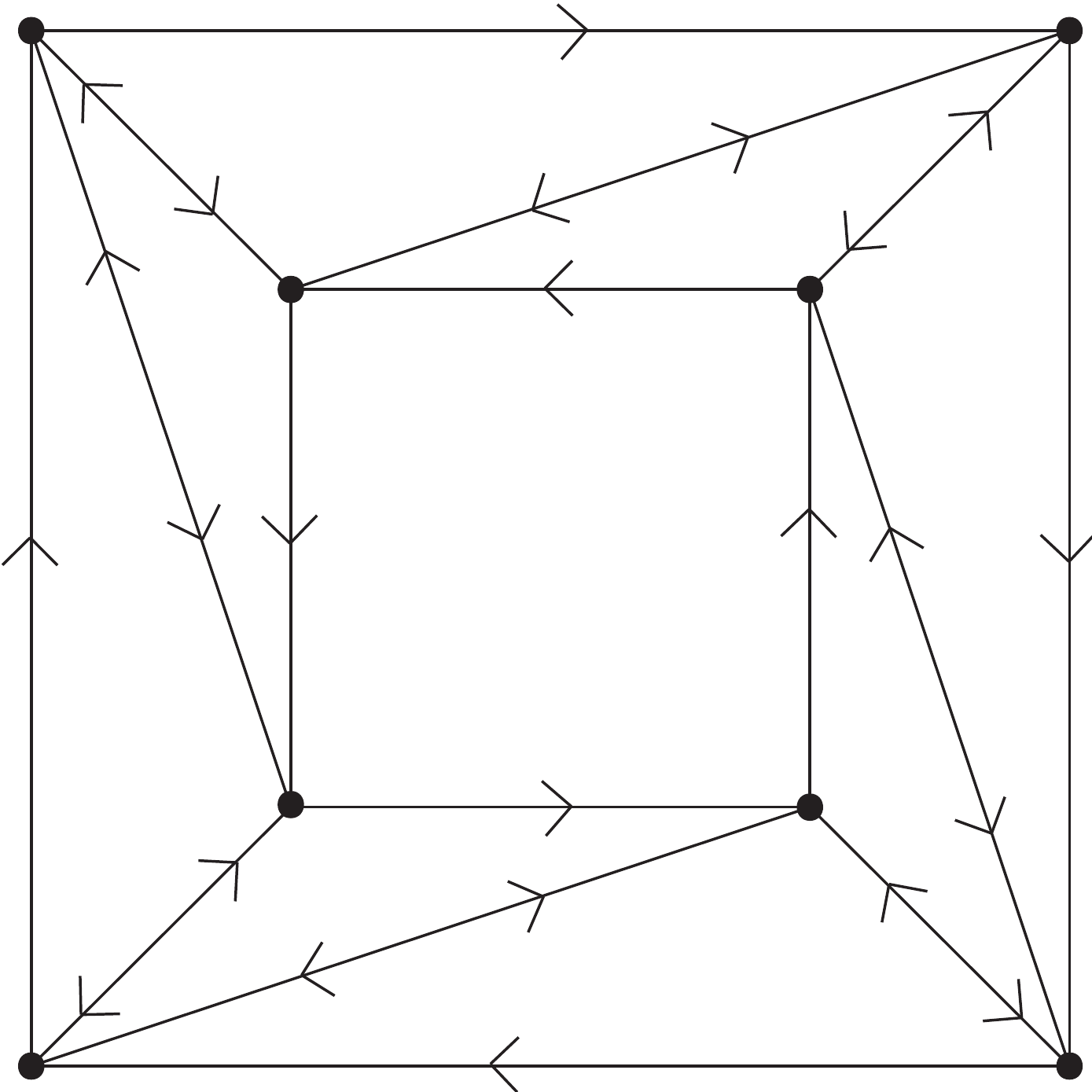}
\center Figure 1
\end{center}


\section{Feasibility Conditions and Known Construction Methods}\label{smain1}

It is immediate from the definition of directed strongly regular
graphs that $0\leq \mu , \lambda$ and $0\leq t\leq k\leq n-1$.
However, we will only consider the (`genuine') DSRGs with $0<t<k$,
excluding totally undirected ($t=k$) and totally directed ($t=0$) cases.
We are able to put further restrictions on parameter sets, which
determines if it is feasible for a DSRG with certain parameters to
exist.

It is known that if $\Gamma$ is a DSRG$(n,k,t,\lambda,\mu)$ with
adjacency matrix $A$, then the complement $\Gamma '$ of $\Gamma$ is
a DSRG$(n,k',t',\lambda ',\mu ')$ with adjacency matrix $A'$ =
$J-I-A$ where $$ k' =(n-2k)+(k-1)$$
$$\lambda ' = (n-2k) + (\mu-2)$$ $$t'=(n-2k)+(t-1)$$ $$\mu
'=(n-2k)+\lambda.$$ This can easily be shown by evaluating $A'^2$.

Duval showed that if a DSRG$(n,k,t,\lambda,\mu)$ is not totally
undirected or complete, nor totally directed, then the following
equations and inequalities hold true:
$$k(k+(\mu-\lambda))=t+(n-1)\mu$$
$$(\mu-\lambda)^2 + 4(t-\mu)=d^2$$
$$d\mid 2k-(\mu-\lambda)(n-1)$$
$$\frac{2k-(\mu-\lambda)(n-1)}{d}
\equiv n-1\ (\mbox{mod } 2)$$ $$\left|
\frac{2k-(\mu-\lambda)(n-1)}{d}\right| \leq n-1$$
$$0\leq\lambda<t<k$$ $$0<\mu\leq t<k$$ $$-2(k-t-1)\leq \mu-\lambda
\leq2(k-t)$$

These parameter restrictions allow for a list of the feasible
directed strongly regular graphs to be compiled.  This list
immediately suggests to the observer that no DSRG of prime order
exists, which is in fact the case (cf. \cite{Du}).

Duval provided an initial list of feasible parameter sets in his
paper \cite{Du}, but a more complete list is available in \cite{BH}.


Duval described many different construction methods including three
that we will describe here:
\begin{enumerate}
\item Constructing directed strongly regular graphs using quadratic residue
\item Constructing DSRGs with a block construction using permutation matrices
\item Constructing DSRGs by using the Kronecker product to construct new DSRGs
from smaller ones
\end{enumerate}

The first construction uses quadratic residue matrices to construct
DSRG$(n,k,t,\lambda,\mu)$ with parameters $(2q, q-1,
\frac{1}{2}(q-1), \frac{1}{2}(q-1) -1, \frac{1}{2}(q-1))$, where
$q=4m+1$ and is a prime power. The adjacency matrices of such DSRGs
will take the form $$A = \left[ \begin{array}{cc}
Q & C_1 \\
C_2 & Q \end{array} \right] .$$

$C_1$ and $C_2$ are $\sigma_1$ and $\sigma_2$ circulant matrices
respectively, where a $\sigma$ circulant matrix $C$ satisfies
$C_{ij}=C_{i-k,j-\sigma{k}}$.  This means that each row, or each
column, is equal to the previous row (column) shifted $\sigma$
entries to the right (down).  $Q$ is a quadratic residue matrix of
order $q$, indexed by the elements of GF($q$), the Galois Field of
order $q$.  When $R$ is the set of quadratic residues of GF($q$),
the nonzero elements $x \in$ GF$(q)$ such that $x=y^2$ for some
$y\in$ GF($q$), and $N$ is the set of quadratic non-residues of
GF($q$), all other nonzero elements of GF($q$), $Q$ is defined by
$$Q_{ij}= \left\{
\begin{array}{ll}
  1 &\mbox{if}\; i - j\; \in R\\
  0 &\mbox{if}\; i - j\; \in N
\end{array}
\right..$$

This construction method produces a DSRG iff
\begin{itemize}
\item $\sigma_1\sigma_2$ = $1 \in GF(q)$
\item $\sigma_1, \sigma_2 \in N$
\item The partition of $GF(q)^*$ into the two sets,
each of $2m$ elements, $$S = \{x \in GF(q)^*:(C_2)_{0,x} =
1\}\quad\mbox{and}\quad T = \{x \in GF(q)^*:(C_2)_{0,x} = 0\}$$
described by the first row satisfies the following ``difference
partition" property: Each of the $4m$ elements of $GF(q)^*$ occurs
exactly $m$ times in the $4m^2$ differences $s-t$ where $s\in S$ and
$t\in T$.
\end{itemize}

An example of an adjacency matrix for the DSRG$(10,4,2,1,2)$ using
the preceding construction is $$A = \left[
\begin{array}{ccccc|ccccc}
0 & 1 & 0 & 0 & 1 & 0 & 1 & 0 & 0 & 1\\
1 & 0 & 1 & 0 & 0 & 0 & 1 & 0 & 1 & 0\\
0 & 1 & 0 & 1 & 0 & 1 & 0 & 0 & 1 & 0\\
0 & 0 & 1 & 0 & 1 & 1 & 0 & 1 & 0 & 0\\
1 & 0 & 0 & 1 & 0 & 0 & 0 & 1 & 0 & 1\\
\hline
0 & 1 & 0 & 0 & 1 & 0 & 1 & 0 & 0 & 1\\
0 & 0 & 1 & 0 & 1 & 1 & 0 & 1 & 0 & 0\\
1 & 0 & 1 & 0 & 0 & 0 & 1 & 0 & 1 & 0\\
1 & 0 & 0 & 1 & 0 & 0 & 0 & 1 & 0 & 1\\
0 & 1 & 0 & 1 & 0 & 1 & 0 & 0 & 1 & 0\\ \end{array} \right] .$$

DSRG's with parameter $(2(2\mu+1), 2\mu, \mu, \mu-1, \mu)$ can be
found by another, simpler, block construction using matrices of the
form $$ A=\left[ \begin{array}{cc}
Q & PQ \\
(PQ)^T & Q \end{array} \right]$$ where $$ Q+Q^T=J-I $$ $$QJ=JQ=\mu
J$$ and $P$ is a permutation matrix with rank 2, so $$PJ=JP=J$$
$$P=P^T=P^{-1}.$$ This construction method yields a DSRG iff $PQ =
(PQ)^T$.

An example of an adjacency matrix for the DSRG$(14,6,3,2,3)$ using
the preceding construction is $$A = \left[\begin{array}{ccccccc|ccccccc}
0 & 1 & 1 & 1 & 0 & 0 & 0 & 0 & 1 & 1 & 1 & 0 & 0 & 0\\
0 & 0 & 1 & 1 & 1 & 0 & 0 & 1 & 1 & 1 & 0 & 0 & 0 & 0\\
0 & 0 & 0 & 1 & 1 & 1 & 0 & 1 & 1 & 0 & 0 & 0 & 0 & 1\\
0 & 0 & 0 & 0 & 1 & 1 & 1 & 1 & 0 & 0 & 0 & 0 & 1 & 1\\
1 & 0 & 0 & 0 & 0 & 1 & 1 & 0 & 0 & 0 & 0 & 1 & 1 & 1 \\
1 & 1 & 0 & 0 & 0 & 0 & 1 & 0 & 0 & 0 & 1 & 1 & 1 & 0\\
1 & 1 & 1 & 0 & 0 & 0 & 0 & 0 & 0 & 1 & 1 & 1 & 0 & 0\\
\hline 0 & 1 & 1 & 1 & 0 & 0 & 0 & 0 & 1 & 1 & 1 & 0 & 0 & 0\\
1 & 1 & 1 & 0 & 0 & 0 & 0 & 0 & 0 & 1 & 1 & 1 & 0 & 0\\
1 & 1 & 0 & 0 & 0 & 0 & 1 & 0 & 0 & 0 & 1 & 1 & 1 & 0\\
1 & 0 & 0 & 0 & 0 & 1 & 1 & 0 & 0 & 0 & 0 & 1 & 1 & 1\\
0 & 0 & 0 & 0 & 1 & 1 & 1 & 1 & 0 & 0 & 0 & 0 & 1 & 1 \\
0 & 0 & 0 & 1 & 1 & 1 & 0 & 1 & 1 & 0 & 0 & 0 & 0 & 1\\
0 & 0 & 1 & 1 & 1 & 0 & 0 & 1 & 1 & 1 & 0 & 0 & 0 & 0 \end{array}\right].$$

A final construction developed by Duval uses the Kronecker product
of matrices.  A small example of how the Kronecker product,
$A\otimes B$, works is shown below. If $$A = \left[
\begin{array}{cc}
a_{11} & a_{12} \\
a_{21} & a_{22} \end{array} \right] $$
and $$B = \left[ \begin{array}{cc}
b_{11} & b_{12} \\
b_{21} & b_{22} \end{array} \right], $$
their Kronecker Product
$$A \otimes B = \left[ \begin{array}{cc}
a_{11}B & a_{12}B \\
a_{21}B & a_{22}B \end{array} \right]=
\left[ \begin{array}{cccc}
a_{11}b_{11} & a_{11}b_{12} & a_{12}b_{11} & a_{12}b_{12} \\
a_{11}b_{21} & a_{11}b_{22} & a_{12}b_{21} & a_{12}b_{22}\\
a_{21}b_{11} & a_{21}b_{12} & a_{22}b_{11} & a_{22}b_{12}\\
a_{21}b_{21} & a_{21}b_{22} & a_{22}b_{21} & a_{22}b_{22}
\end{array} \right].$$
The construction method works as follows: let $A$ be the adjacency
matrix of a DSRG and $J_m$ be the all-ones matrix.  For $m>1$, $A \otimes J_{m}$ is the adjacency
matrix of a DSRG$(nm, km, tm, \lambda{m},  \mu{m})$ iff $t=\mu$.  a
DSRG with the same parameters as above can also be constructed from
$J_{m} \otimes A$.

\vspace{.1 in}We will discuss some other construction methods not introduced by Duval.  We begin with an  introduction to constructing DSRGs using Cayley graphs. Let $G$ be a finite group and $S \subset G-\{e\}$.  The Cayley Graph of $G$ generated by $S$, $Cay(G; S)$, is the digraph $\Gm$ such that $V(\Gamma)=G$ and \[E(\Gamma) = \{(x,y): \exists \;s\in S \mbox{ such that } xs=y\}.\]

An example of a Cayley graph that is from an abelian group is illustrated in Figure 2. An example of a
DSRG constructed from a Cayley graph is illustrated in Figure 3.
\begin{center}
\includegraphics[scale=0.3]{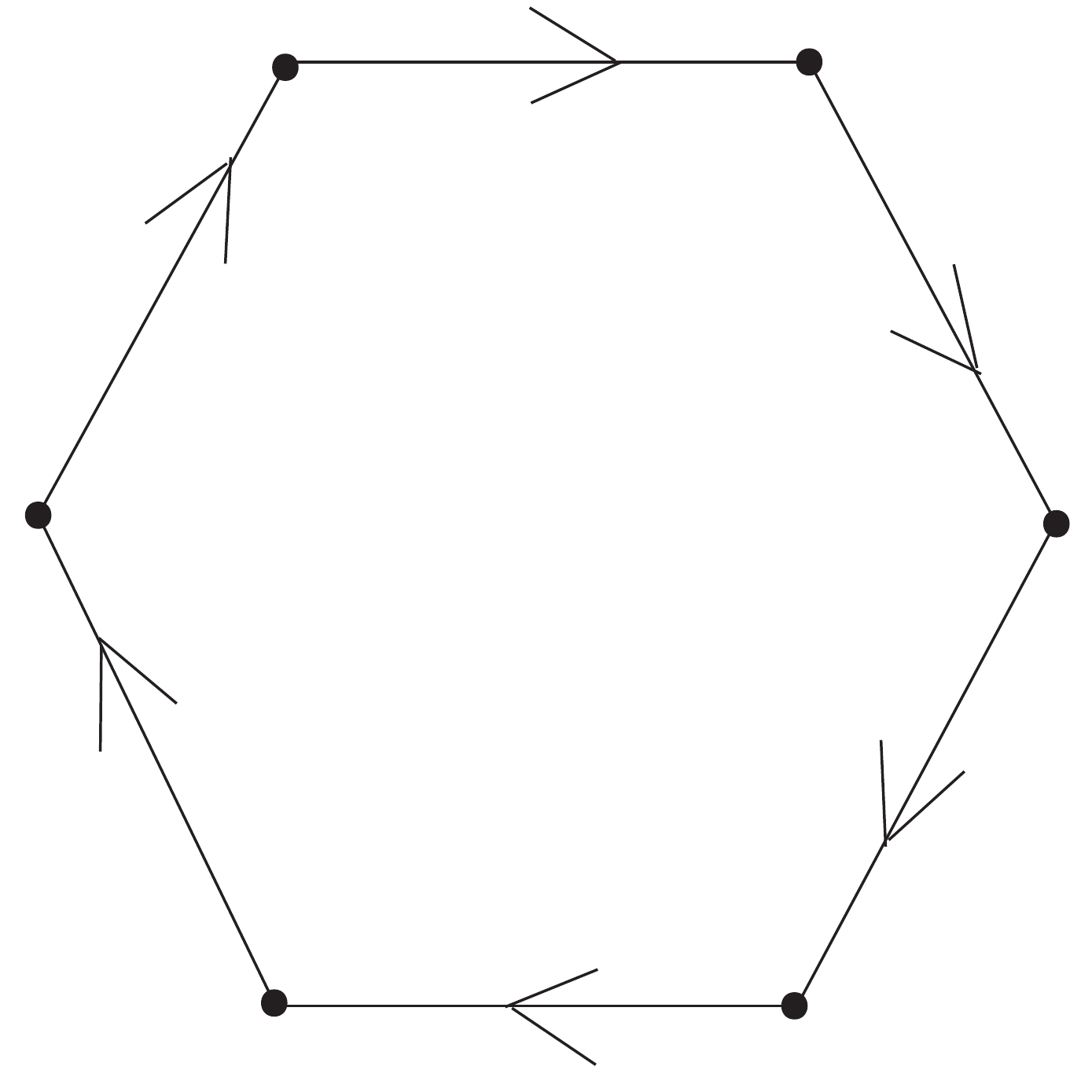}
\center Figure 2, Cay$(\Z_6, \{1\})$
\end{center}

\begin{center}
\includegraphics[scale=0.25]{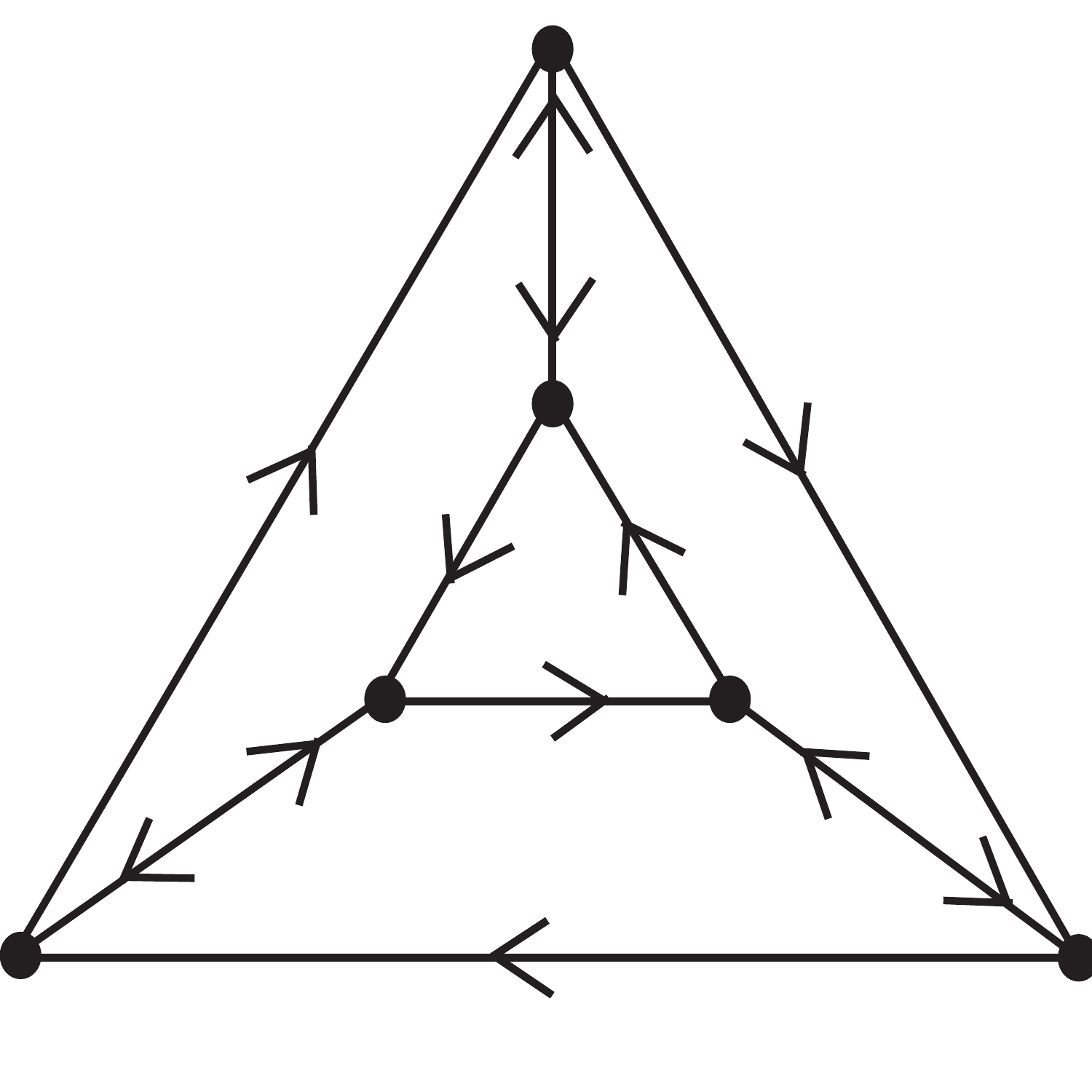}
\center Figure 3, Cay$(S_3, \{(12),\;(123)\})$
\end{center}

Cayley graphs of groups give us another way to construct DSRGs. It
is clear that if Cay$(G; S)$ is to be a DSRG$(n, k, t, \lambda,
\mu)$, $|G| = n$ and $|S| = k$.   Also, it is necessary that in the
induced multiplication table of $S$:
\begin{itemize}
\item The identity of the group, $e$, appears $t$ times.
An easier way to check this is, if $S^{-1}$ = \{$x\in G : x^{-1} \in
S$\}, then $|S \cap S^{-1}| = t$
\item Each element of $S$ appears $\lambda$ times.
\item Each of the elements of $G-S-\{e\}$ appears $\mu$ times.
\end{itemize}
 A result by J{\o}rgensen shows that if $G$ is abelian,
 then Cay$(G; S)$ is not a DSRG for any $S\subset G$ \cite{Jo}.

Other constructions of DSRGs as Cayley graphs were developed by Hobart
and Shaw \cite{HS}. They used the dihedral group $D_{2n}=\;<\alpha,
\beta:\beta^2=\alpha^n=e \mbox{ and }\beta\alpha\beta=\alpha^{-1}>$.
They showed how Cay$(D_{2n}, S)$ can be a DSRG as follows:
\begin{enumerate}
\item When $n=2\lambda$, an even integer, they constructed the DSRG
$(4\lambda, 2\lambda -1, \lambda, \lambda-1, \lambda-1)$  by setting
$S=\{\alpha, \alpha^2,\dots, \alpha^{\lambda-1}, \beta, \beta\alpha,
\dots, \beta\alpha^{\lambda-1}\}$.
\item For $n=2\lambda +1$, an odd integer, they constructed a
DSRG$(4\lambda +2, 2\lambda+1, \lambda, \lambda-1, \lambda)$ by
letting $S=\{\alpha, \alpha^2,\dots, \alpha^{\lambda}, \beta,
\beta\alpha,\dots, \beta\alpha^{\lambda}\}$.
\end{enumerate}
These graphs also appear in \cite{KM}.  
\vspace{.1 in}  
Analyzing the adjacency matrices of these graphs , we observe that these adjacency matrices can be expressed as block matrices of the form
\[B = \left [
\begin{array}{cc}
A & A^T \\
A & A^T \end{array} \right ] \] where $A$ is the adjacency matrix of
a highly structured graph, namely, a regular tournament. We observe
that this form of block matrix is indeed able to be used in a
general construction method for DSRGs.

\section{New Constructions using Tournaments and Circulant Matrices }\label{s3}

In this section we will introduce several new methods of
constructing directed strongly regular graphs that fall into three
categories.
\begin{enumerate}
\item Constructing DSRGs using regular tournaments
\item Constructing DSRGs using doubly regular tournaments
\item Constructing DSRGs using circulant matrices
\end{enumerate}


\begin{defn}\label{defn1} A tournament is a directed graph $\Gamma$
such that for any $x,y \in V(\Gamma)$ exactly one of $x\rightarrow
y$ or $y\rightarrow x$ holds.  A tournament $\Gm$ is said to be
regular if every vertex in $V(\Gm)$ has the same out-degree. Thus a
regular tournament has $n=2k+1$ if $n$ and $k$ denote the number of
vertices and the valency of the graph, respectively.
\end{defn}

The adjacency matrix $A$ of a tournament $\Gm$ satisfies the equation
$A+A^T=J-I$.  If $\Gm$ is a regular tournament with valency $k$,
then $JA=AJ=kJ$.


\begin{lem}\label{lem2} If $A$ is an adjacency matrix of a regular
tournament with valency $k$, then

\begin{enumerate}
\item $B = \left[ \begin{array}{cc}
A & A^T \\
A & A^T \end{array} \right] $
\item $C =
\left[\begin{array}{cc}
A & A \\
A^T & A^T \end{array} \right] $
\end{enumerate}
are adjacency matrices of directed
strongly regular graphs with parameters $(4k+2, 2k, k, k-1, k)$.
\end{lem}

\bpf Let $J$ denote the  $(4k+2) \times (4k+2)$ all-ones matrix
while $\bar{J}$ denote the $(2k+1) \times (2k+1)$ all-ones matrix,
and similarly for $I$ and $\bar{I}$. Then we have
$$JB = BJ = \left[ \begin{array}{cc}
A & A^T \\
A & A^T \end{array} \right)
\left( \begin{array}{cc}
\bar{J} & \bar{J} \\
\bar{J}& \bar{J} \end{array} \right] = 2kJ.$$ Since
$A\bar{J}=\bar{J}A=k\bar{J}$ and $A + A^T + \bar{I}=\bar{J}$,
$$B^2 + B=\left[\begin{array}{c|c}
{A^{2} + A^{T}A + A} & {AA^{T} + (A^T)^{2} + A^{T}}\\
\hline{A^{2} + A^{T}A + A} & {AA^{T} + (A^T)^{2} + A^{T}}
\end{array} \right]=kJ.$$ Equivalently, $$B^2=kI + (k-1)B +
k(J-I-B).$$ Therefore $B$ is an adjacency matrix of a DSRG$(4k+2,
2k, k, k-1, k)$.  Similarly, it can be shown that matrix $C$ is also
the adjacency matrix of a DSRG$(4k+2, 2k, k, k-1, k)$. \epf

A similar block construction also produces DSRGs and is actually closely related to the previous construction.

\begin{lem}\label{lem8} If the matrix $A$ is the adjacency matrix of
a regular tournament of order $2k+1$, then the matrix \[M(A)=
\left[\begin{array}{cc}
A & A^T+I\\
A+I & A^T \end{array} \right]\] is the adjacency matrix of a DSRG
with parameters $(4k+2, 2k+1, k+1, k, k+1)$.
\end{lem}

\bpf From Lemma \ref{lem2}, we know that the adjacency matrix $B =
\left[ \begin{array}{cc}
A & A^T\\
A & A^T \end{array} \right]$ is the adjacency matrix of a DRSG with
parameters $(4k+2, 2k, k, k-1, k)$.  Duval showed that the
complement of a DSRG is also a DSRG (see Section \ref{smain1}).  We
will show that graph $\Gm'$ represented by $M$ is the complement of
the graph $\Gm$ represented by $B$.

Since $A$ is an adjacency matrix of a regular tournament, it holds
$A+A^T={J}-{I}$. Therefore, the complement $B'$ of $B$ can be
simplified to

$$\left[ \begin{array}{cc}
J-I-A & J-A^T\\
J-A & J-I-A^T \end{array} \right]  = \left[
\begin{array}{cc}
A^T & A+{I}\\
A^T +{I} & A \end{array} \right].$$ If we choose $P$ to be the
permutation matrix equal to $\left[\begin{array}{cc}
\bz & {I}\\
{I} & \bz \end{array} \right] $  where $\bz$ denotes the all-zeros
matrix. Then
$$PB'P = \left[\begin{array}{cc}
A & A^T +I \\
A+I & A^T \end{array} \right]= M.$$ Therefore, $B'\cong M$ and $M$
is the adjacency matrix of DSRG with parameters $(4k+2, 2k+1, k+1,
k, k+1)$.  A completely analogous construction coming from Lemma \ref{lem2}.2 also produces DSRGs in a nice block matrix form.

\epf
\vspace{.1 in} To simplify our next few constructions, we will use the notation $M(A)$ to mean the matrix of the form 
$\left[\begin{array}{cc}
A & A^T+I\\
A+I & A^T \end{array} \right]$ for any matrix $A$.


This construction method can be generalized by having multiple
columns or rows of $A$ and $A^T$ for a regular tournament.

\begin{lem}\label{lem4} If $A$ is an adjacency matrix of a regular
tournament with valency $k$, then
\begin{enumerate}
  \item $B = \left[\begin{array}{ccccc}
A & A^T & A & \ldots & A^T \\
A & A^T & A & \ldots & A^T \\
\vdots & \vdots & \vdots & \ddots & \vdots \\
A & A^T & A & \ldots & A^T  \end{array} \right] $
  \item $ C = \left[ \begin{array}{cccc}
A & A & \ldots & A \\
A^T & A^T & \ldots & A^T \\
A & A &\ldots & A \\
\vdots & \vdots  & \ddots & \vdots \\
A^T & A^T & \ldots & A^T   \end{array} \right]$
\end{enumerate}
are adjacency matrices of DSRGs with parameters $((4k +2)w, 2kw, kw,
(k-1)w, kw)$ where $w$ is the number of $A$ and $A^T$ blocks in each
column (or row).
\end{lem}

\bpf
Since $B = \left[ \begin{array}{cc}
A & A^T \\
A & A^T \end{array} \right] $ is the adjacency matrix of a
DSRG$(4k+2, 2k, k, k-1, k)$, using Duval's Kronecker product
construction, $J_w\otimes B$ is also the adjacency matrix of a DSRG,
but with parameters $((4k +2)w, 2kw, kw, (k-1)w, kw)$.  The proof
that $C$ is a DSRG$(4k+2, 2k, k, k-1, k)$ is exactly the same. \epf


\begin{defn}\label{defn2}  A regular tournament $T$
is said to be \textit{doubly regular} if for every vertex $x\in
V(T)$, the out-neighbors of $x$ span a regular tournament. If $T$ is
a doubly regular tournament of order $n$, with regular valency $k$
and the degree of the induced subgraph on the out-neighbors
$\lambda$, then $n=2k+1=4\lambda+3$. 
\end{defn}


\begin{defn}\label{defn3} An $(m,r)$-team tournament is a digraph obtained from
the complement $\overline{m\circ K_r}$ of $m$ copies of the
complete graph $K_r$ by giving an orientation in such a way that
every undirected edge $\{x,y\}$ is assigned with either
$x\rightarrow y$ or $x\leftarrow y$ but not both. 
\end{defn}

We note that an $(m,r)$-team tournament has $m$ maximal independent
sets of size $r$, and the edges are directed links between the
vertices of distinct maximal independent sets.

\begin{defn}\label{defn4} An $(m,r)$-team tournament $\Gm$
with adjacency matrix $A$ is said to be \textit{doubly regular} iff
\bit
\item[(1)] every vertex of $\Gm$ has in-degree and out-degree
$k=\frac{1}{2}(m-1)r$, and
\item[(2)] there exist positive integers $\alpha, \beta$ and $\gamma$ such that
for every pair of distinct vertices $x$ and $y$, the number of
directed paths of length $2$ from $x$ to $y$ is
$$
\left\{
\begin{array}{ll}
  \alpha & \mbox{if }\; x \rightarrow y\\
 \beta & \mbox{if }\; x \leftarrow y\\
 \gamma & otherwise
\end{array}
\right ..$$ \eit
\end{defn}
As in \cite{JJ}, the adjacency matrix of a doubly regular
$(m,r)$-team tournament satisfies the following equations\bit
\item[(1)] $AJ=JA=kJ$;
\item[(2)] $A^2 = \alpha{A} + \beta{A^T} +\gamma(J-I-A-A^T)$.
\eit

In \cite{JJ}, we can find $(m,2)$-team tournaments coming from
doubly regular tournaments of order $m-1$. Let $A$ be an adjacency
matrix of a doubly regular tournament $T$ of order
$m-1=2k+1=4\lambda+3$. Then
$$D(T) =  \left[ \begin{array}{c|ccc|c|ccc}
0 & 1 & \ldots & 1& 0 & 0 & \ldots & 0 \\
\hline 0 &    &     &   & 1 &     &     &    \\
:   &    & A &   & :  &     & A^T& \\
0 &    &     &   & 1 &     &     &    \\
\hline 0 & 0 & \ldots & 0& 0 & 1 & \ldots & 1 \\
\hline 1 &    &     &   & 0 &     &     &    \\
:   &    & A^T &   & :  &     & A& \\
1 &    &     &   & 0 &     &     &    \\
 \end{array} \right] $$
is an adjacency matrix of a doubly regular $(m, 2)$-team
tournament.


\begin{lem}\label{lem5}  Let $D=D(T)$ be an $(m,2)$-team tournament
described above with $m=2k+2=4\lambda +4$, then $$M =M(D)= \left[
\begin{array}{cc}
D & D^T+I \\
D+I & D^T \end{array} \right]$$ is an adjacency matrix of a DSRG
with parameters $$(4m, 2m-1, m, m-1, m-1)=(16\lambda +16,
8\lambda+7, 4\lambda+4, 4\lambda+3, 4\lambda+3).$$
\end{lem}

\bpf Being an adjacency matrix of doubly regular $(m, 2)$-team
tournament, it is shown in \cite{JJ} that $D$ satisfies \bit
\item[(a)] $D^2=(2\lambda+1)(D+D^T)+(4\lambda+3)(J-I-D-D^T)$,
\item[(b)] $DD^T=(4\lambda+3)I +(2\lambda +1)(D +D^T)$.
\eit 
First we will show that $MJ=JM=kJ=(2m-1)J$
$$MJ=\left[\begin{array}{cc}
D\bar{J}+D^{T}\bar{J} + \bar{J} & D\bar{J}+D^{T}\bar{J} + \bar{J} \\
D\bar{J}+D^{T}\bar{J} + \bar{J} & D\bar{J}+D^{T}\bar{J} + \bar{J} \end{array}\right]=JM$$

Using the definition of a regular tournament we can see that since the doubly regular tournament D has valency $m-1$
$$MJ=JM=\left[\begin{array}{cc}
\bar{J}(m-1+m-1+1) & \bar{J}(m-1+m-1+1)  \\
\bar{J} (m-1+m-1+1) & \bar{J}(m-1+m-1+1)  \end{array}\right]= (2m-1)J$$

Therefore the equation $MJ=JM=kJ$ holds.

We now square the adjacency matrix $M$ to show that $M$ is an
adjacency matrix of the desired DSRG.
$$M^2= \left[ \begin{array}{c|c}
D^2 + D^{T}D + D + D^T + I & DD^T + D + D^T + {D^T}^2\\
\hline D^2 + D^{T}D + D + D^T & DD^T + D + D^T + {D^T}^2 + I
\end{array} \right].$$ Simplifying this by using the equalities
(a) and (b) as well as $D+D^T=J-I$ yields $$M^2 = \left[
\begin{array}{c|c}
(4\lambda +3)J + I & (4\lambda +3)J \\
\hline (4\lambda +3)J & (4\lambda +3)J + I \end{array} \right].$$
Therefore, $t=4\lambda +4, \lambda=\mu=4\lambda +3$.  Each vertex
will have in- and out-valency equal to $2m-1$ because each doubly
regular $(m,2)$-team tournament has valency $m-1$.  Each DSRG will
have $4m$ vertices because in a doubly regular $(m,2)$-team
tournament the number of vertices equals $2m$ and in the adjacency
matrix $M$, thus $v=4m$.  Therefore, we have shown that we can use a
doubly regular $(m, 2)$ team tournament to construct a DSRG with the
parameters $(4m, 2m-1, m, m-1, m-1)$ where $m\equiv 0 \; (\mbox{mod}
\;4)$. \epf


We can extend this construction in lemma \ref{lem5} to make use of
any regular tournament instead of only doubly regular tournaments.

\begin{lem}\label{lem6}
Let $A$ be the adjacency matrix of a regular tournament $T$ of order
$h$ and
$$D=D(T) = \left[ \begin{array}{cccc}
0 & \boT & 0 & \bzT \\
\bz & A & \bo& A^T \\
0 & \bzT & 0 & \boT\\
\bo& A^T &  \bz & A\\
\end{array} \right] $$
where $\bz$ and $\bo$ denote the $n$-dimensional column vectors of
all zeros and all ones, respectively. The matrix $$M(D) = \left[
\begin{array}{cc}
D & D^T+I \\
D+I & D^T \end{array} \right]$$ is the adjacency matrix of a DSRG
$(4(h+1), 2h+1, h+1, h, h)$ where $h\equiv 1 \mbox{ mod }2$.
 \end{lem}

\bpf 
First we will show that $MJ=JM=kJ=(2h+1)J$
$$MJ=\left[\begin{array}{cc}
D\bar{J}+D^{T}\bar{J} + \bar{J} & D\bar{J}+D^{T}\bar{J} + \bar{J} \\
D\bar{J}+D^{T}\bar{J} + \bar{J} & D\bar{J}+D^{T}\bar{J} + \bar{J} \end{array}\right]=JM$$

Using the definition of a regular tournament we can see that since the regular tournament D has valency $m$ then 
$$MJ=JM=\left[\begin{array}{cc}
\bar{J}(h+h+1) & \bar{J}(h+h+1)  \\
\bar{J} (h+h+1) & \bar{J}(h+h+1)  \end{array}\right]= (2h+1)J$$

Therefore the equation $MJ=JM=kJ$ holds.

Since $A+A^T=J_h-I_h$, we have
$$D +D^T\ = \left[ \begin{array}{cc|cc}
0 & \boT & 0 & \boT \\
\bo & A + A^T & \bo& A+A^T \\
\hline 0 & \boT & 0 & \boT\\
\bo& A+A^T &  \bo & A+A^T\\
\end{array} \right]\ =\ \left[ \begin{array}{c|c}
J_{h+1}-I_{h+1} & J_{h+1}-I_{h+1}\\
\hline J_{h+1}-I_{h+1} & J_{h+1}-I_{h+1} \end{array} \right].$$

$$D^{T}D= DD^T=  \left[ \begin{array}{c|c|c|c}
h & k\boT & 0 & k\boT \\
\hline k\bo & J_{h}+A^{T}A+AA^T & k\bo& A^2 + {A^T}^2 \\
\hline 0 & k\boT & n & k\bo\\
\hline k\bo& A^2+{A^T}^2 &  k\bo & J_{n} + AA^T+A^{T}A\\
 \end{array} \right] $$

\[D^2 = (D^T)^2 =  \left[ \begin{array}{c|c|c|c}
0 & k\boT & h & k\boT \\
\hline k\bo & A^2 + {A^T}^2 & k\bo& J_{h}+A^{T}A+AA^T \\
\hline h & k\boT & 0 & k\bo\\
\hline k\bo& J_{h}+A^{T}A+AA^T &  k\bo & A^2+{A^T}^2 \\
 \end{array} \right] \]

$D^{2} +DD^{T} + D +D^{T}=$ $$ \left[ \begin{array}{c|c|c|c}
h & (2k+1)\boT & h & (2k+1)\boT \\
\hline (2k+1)\bo & (A+A^T)^2 + J_{h}+ A+A^T & n\bo& (A+A^T)^2 + J_{h}+ +A +A^T \\
\hline h & (2k+1)\boT & h & (2k+1)\bo\\
\hline (2k+1)\bo& (A+A^T)^2 + J_{h}+ A +A^T & h\bo & (A+A^T)^2 + J_{h}+ A + A^T \\
 \end{array} \right] $$

Using the fact that for a regular tournament of order $h$ with
adjacency matrix $A$, $A + A^T=J_h-I_h$ we can easily
simplify $D^{2}+D^{T}+ D +D^{T}$ to $hJ_{2h+2}$

From  lemma \ref{lem5}, it is known that
$$M(D)^2= \left[ \begin{array}{c|c}
D^2 + D^{T}D + D + D^T + I & DD^T + D + D^T + {D^T}^2\\
\hline D^2 + D^{T}D + D + D^T & DD^T + D + D^T + {D^T}^2 + I
\end{array} \right] .$$

Using the above simplification we can transform
$$M^2= \left[ \begin{array}{c|c}
hJ_{2h+2} +I_{2h+2} & hJ_{2h+2}\\
\hline hJ_{2h+2} & hJ_{2h+2} + I_{2h+2} \end{array} \right]$$
into
$$hJ_{4n+4} +I_{4h+4}$$

Therefore $M$ is the adjacency matrix of the DSRG with parameters
$(4(h+1), 2h+1, h+1, h, h)=(8(k+1), 4k+3, 2k+2, 2k+1, 2k+1)$ for
$k=\frac12(h-1)\in \Z^+$. \epf

In what follows, let $\Pi=\Pi(n)$ denote an $n\times n$ permutation
matrix corresponding to the $n$-cycle $(1, 2 \cdots n)$ given by
$$\Pi= \left[ \begin{array}{cccccc}
0 & 1 & 0 & 0 & \ldots & 0 \\
0 & 0 & 1 & 0 & \ldots & 0 \\
0 & 0 & 0 & 1 & \ldots & 0 \\
\vdots & \vdots & \vdots & \ddots & \ddots & \vdots\\
\vdots & \vdots & \vdots &  &  \ddots & 1\\
1 & 0 & 0 & \ldots & \ldots & 0
 \end{array} \right].$$ Then $\Pi^2$ is the permutation matrix of
another $n$-cycle, but with the diagonal moved up to the second
off-diagonal, and so on, until $\Pi^{n}=I$.

\begin{lem}\label{lem7} For a positive integer $s$, let $L$ be the
$(2s+2)\times (2s+2)$-matrix equals to $\Pi + \Pi^2 + \cdots +
\Pi^s$ where $\Pi=\Pi(2s+2)$ as defined above. Then $$M(L) = \left[
\begin{array}{cc}
L & L^T+I \\
L+I & L^T \end{array} \right] $$ is an adjacency matrix of a DSRG
with parameters $(4(s+1), 2s+1, s+1, s, s)$.
\end{lem}

\bpf 

First we will show that $MJ=JM=kJ=(2s+1)J$
$$MJ=\left[\begin{array}{cc}
L\bar{J}+L^{T}\bar{J} + \bar{J} & L\bar{J}+L^{T}\bar{J} + \bar{J} \\
L\bar{J}+L^{T}\bar{J} + \bar{J} & L\bar{J}+L^{T}\bar{J} + \bar{J} \end{array}\right]=JM$$

Since we know that the matrix L has s ones per row we can transform 

$$MJ=JM=\left[\begin{array}{cc}
\bar{J}(s+s+1) & \bar{J}(s+s+1)  \\
\bar{J} (s+s+1) & \bar{J}(s+s+1)  \end{array}\right]= (2s+1)J$$

Therefore the equation $MJ=JM=kJ$ holds.

The matrix $L$ may be expressed as the block matrix $$L=\left[
\begin{array}{cc}
B & B^T\\
B^T & B \end{array} \right] $$ with the $(s+1) \times (s+1)$-matrix
$$B = \left[ \begin{array}{ccccc}
0 & 1 & \ldots & \ldots & 1\\
0 & 0 & 1 & \ldots & 1\\
\vdots & & \ddots &  & \vdots\\
\vdots & & & \ddots & 1 \\
0 & \cdots & \cdots  & \cdots & 0
 \end{array} \right] $$ the matrix with an all one upper triangle.
It is easy to see that
$$L+L^T = \left[ \begin{array}{c|c}
B + B^T & B +B^T\\
\hline B +B^T & B + B^T \end{array} \right] =  \left[
\begin{array}{c|c}
J-I & J-I\\ \hline
J-I & J-I \end{array} \right]; $$

$$L^2 = (L^T)^2 = \left[ \begin{array}{cc}
B^2 +{B^T}^2 & BB^T + {B^T}B \\
BB^T + {B^T}B  & B^2 +{B^T}^2  \end{array} \right]; $$

$$LL^T= L^{T}L = \left[ \begin{array}{cc}
BB^T + {B^T}B & B^2 +{B^T}^2\\
B^2 +{B^T}^2 & BB^T + {B^T}B\end{array} \right]; $$ and thus,
$$L^{2} +LL^{T} + L +L^{T} = sJ_{2s+2} $$
As in Lemma \ref{lem5},
$$M^2= \left[ \begin{array}{c|c}
L^2 + L^{T}L + L + L^T + I & LL^T + L + L^T + {L^T}^2\\
\hline L^2 + L^{T}L + L + L^T & LL^T + L + L^T + {L^T}^2 + I
\end{array} \right].$$
Using the above simplification, we can transform $M^2$ into
$$M^2= \left[ \begin{array}{c|c}
sJ_{2s+2} +I_{2s+2} & sJ_{2s+2}\\
\hline sJ_{2s+2} & sJ_{2s+2} + I_{2s+2} \end{array}
\right]=sJ_{4s+4} +I_{4s+4}.$$ Therefore, $M$ is the adjacency
matrix of a DSRG with parameters $(4(s+1), 2s+1, s+1, s, s)$. \epf

\section{Isomorphisms between Constructions}\label{s6}
We have discussed many different construction methods in searching
for new DSRGs. We have seen that graphs with the same parameters
can be obtained from different construction methods. In this section
we investigate whether our construction methods are well-defined in the sense of isomorphic tournaments producing isomorphic graphs and
whether the graphs having the same parameters are isomorphic.

\begin{thm}\label{thm1}
All constructions of DSRGs discussed in Lemmas \ref{lem2},
\ref{lem8}, \ref{lem4},\ref{lem5}, \ref{lem6} and \ref{lem7} are
well defined constructions.
\end{thm}

\bpf
\begin{itemize}
\item[] Lemmas \ref{lem2} and \ref{lem8}:
If $A\cong B$; that is, if there exists a permutation matrix $P$
such that $PAP^{-1}=B$, then
$$\left[
\begin{array}{cc}
P & 0 \\
0 & P \end{array} \right] \left[ \begin{array}{cc}
A & A^T \\
A & A^T \end{array}\right]\left[\begin{array}{cc}
P ^{-1}& 0 \\
0 & P^{-1} \end{array} \right)=\left[ \begin{array}{cc}
B & B^T \\
B & B^T \end{array} \right] $$ so $\left[\begin{array}{cc}
A & A^T \\
A & A^T \end{array} \right]\cong \left[ \begin{array}{cc}
B & B^T \\
B & B^T \end{array} \right],$ and similarly, for Lemma \ref{lem2}.  The same permutation
matrix works for Lemma \ref{lem8}.  Since $$\left[
\begin{array}{cc}
P & 0 \\
0 & P \end{array} \right) \left[ \begin{array}{cc}
A & A^T+I \\
A+I & A^T \end{array}\right] \left[ \begin{array}{cc}
P ^{-1}& 0 \\
0 & P^{-1} \end{array} \right]=\left[ \begin{array}{cc}
B & B^T+I \\
B+I & B^T \end{array} \right] $$  the constructions again produce
isomorphic graphs.
This also shows that for any matrices $A,B$ where $A\cong B$, $M(A)\cong M(B)$.

\item[] Lemma \ref{lem4}:
If $PAP^{-1}=B$, the obvious generalization of the smaller case,
$$\left[\begin{array}{ccccc}
P & 0 & 0 & \ldots & 0 \\
0 & P & 0 & \ldots & 0 \\
0 & 0 & P & \ldots & 0 \\
\vdots & \vdots & \vdots & \ddots & \vdots \\
0 & 0 & 0 & \ldots & P  \end{array} \right]
\left[ \begin{array}{ccccc}
A & A^T & A & \ldots & A^T \\
A & A^T & A & \ldots & A^T \\
\vdots & \vdots & \vdots & \ddots & \vdots \\
A & A^T & A & \ldots & A^T  \end{array} \right]
\left[\begin{array}{ccccc}
P^{-1} & 0 & 0 & \ldots & 0 \\
0 & P^{-1} & 0 & \ldots & 0 \\
0 & 0 & P^{-1} & \ldots & 0 \\
\vdots & \vdots & \vdots & \ddots & \vdots \\
0 & 0 & 0 & \ldots & P^{-1}  \end{array} \right]$$
$$=\left[\begin{array}{ccccc}
B & B^T & B & \ldots & B^T \\
B & B^T & B & \ldots & B^T \\
\vdots & \vdots & \vdots & \ddots & \vdots \\
B & B^T & B & \ldots & B^T  \end{array} \right]$$
shows that again we have isomorphic graphs.

\item[] Lemmas \ref{lem5}, \ref{lem6} and \ref{lem7}:
If we show that given $A\cong B$, $D(A)\cong D(B)$ and $M(A)\cong
M(B)$, it follows that these constructions produce isomorphic graphs
given isomorphic tournaments, or block matrices in the case of Lemma
\ref{lem7}.  Since the proof of Lemma \ref{lem8} shows $M(A)\cong
M(B)$, it only remains to show $D(A)\cong D(B)$. This is done by
seeing that
$$\left[ \begin{array}{c|c|c|c}
1 &  &  &  \\
\hline & P & &  \\
\hline  &  & 1 & \\
\hline &  &  & P\\
 \end{array} \right]
 \left[ \begin{array}{c|c|c|c}
0 & \boT & 0 & \bzT \\
\hline \bz & A & \bo& A^T \\
\hline 0 & \bzT & 0 & \boT\\
\hline \bo& A^T &  \bz & A\\
 \end{array} \right]
 \left[ \begin{array}{c|c|c|c}
 1&  &  &  \\
\hline & P^{-1} & & \\
\hline  &  & 1 & \\
\hline &  &  & P^{-1}\\
 \end{array} \right] $$
 $$ = \left[ \begin{array}{c|c|c|c}
0 & \boT & 0 & \bzT \\
\hline \bz & B & \bo& B^T \\
\hline 0 & \bzT & 0 & \boT\\
\hline \bo& B^T &  \bz & B\\
 \end{array} \right]. $$
\end{itemize}
\epf So, using the same construction, isomorphic tournaments always
produce isomorphic DSRGs, but it is not as clear if, for one
construction, non-isomorphic tournaments always yield non-isomorphic
DSRGs.

With the two very similar constructions from Lemma \ref{lem2}, it is
natural to wonder when, given the same tournament, these two
different constructions actually produce isomorphic graphs.  Below
is one criterion for determining if they are isomorphic.

\begin{lem}\label{lem3} Let $A$ be an adjacency matrix of a regular
tournament and $B, C$ be the adjacency matrices of the DSRGs
constructed from Lemmas \ref{lem2}(1) and
\ref{lem2}(2), respectively.  If there exists a permutation matrix
$P$ such that $P{A}=A^T=AP$, then $B\cong C$.

\end{lem}

\bpf
Assuming $P{A}=A^T=AP$, so $P^{-1}{A^T}=A=A^TP^{-1}$, if
$$\left[ \begin{array}{cc}
I  & 0 \\
0 & P \\
 \end{array} \right)
 \left( \begin{array}{cc}
A & A^T \\
A & A^T \end{array}\right]
\left[\begin{array}{cc}
I  & 0 \\
0 & P^{-1} \\
 \end{array} \right]
= \left[\begin{array}{cc}
A  & A^TP^{-1} \\
P{A} & P{A}P^{-1} \\
 \end{array} \right]=
 \left[ \begin{array}{cc}
A & A \\
A^T & A^T \end{array} \right] = C. $$

So $B\cong C$
\epf

Since circulant matrices commute with each other, this condition is
reduced to $P{A}=A^T$, which happens iff each row of $A$ appears
also as a column of $A$.  This property shows up in many tournaments
that can be decomposed as a sum of permutation matrices representing
$h$-cycles where $|A|=2h+1$.

\begin{prop}\label{prop1} Let $\Pi=\Pi(2k+1)$.
 The matrices defined as $P=\Pi+\Pi^3+\Pi^5+\cdots+\Pi^{2k-1},
 P_0=\Pi^2+\Pi^4+\Pi^6+\cdots+\Pi^{2k},$ and $
 P_j=\Pi^j(\Pi+\Pi^2+\Pi^3+\ldots+\Pi^{k})$ for any $0\leq j\leq k$.
\end{prop}

Just as Lemma \ref{lem4} generalizes our construction in Lemma
\ref{lem2}, the sufficient condition for an isomorphism between the
two DSRGs constructed from a regular tournament $A$ in Lemma
\ref{lem3} holds in the more general case.

\begin{lem}\label{lem9} Let $A$ be an adjacency matrix of a regular
tournament and $B, C$ be the adjacency matrices of the DSRGs
constructed from Lemmas \ref{lem5}(1) and
\ref{lem5}(2), respectively.  If there exists a permutation matrix
$P$ such that $P{A}=A^T=AP$, then $B\cong C$.

\end{lem}

\bpf
 If $$H = \left[\begin{array}{ccccc}
I & 0 & 0 & \ldots&0 \\
0 & P & 0& \ldots & 0 \\
0 & 0 & I &\ldots&0\\
\vdots & \vdots & \vdots & \ddots & \vdots \\
0 & 0 & 0 & \ldots& P \end{array} \right] $$ where $0$ is the matrix
with all-$0$ entries, then $HBH^{-1} = C$, since each 2-block
$\times$ 2-block section reduces to the matrix from Lemma
\ref{lem3}. \epf

It is interesting to note that the block columns of $B$ (block rows
of $C$) could actually be arranged in any order and would still
produce a DSRG isomorphic to the original graph constructed in Lemma
\ref{lem5}.  This isomorphism is associated with the permutation
matrix $P\otimes I_{2h+1}$, for some permutation matrix P of order
$2w$.

\section{Areas of Further Investigation and Summary}\label{s5}

We conclude the paper by making a few remarks.
\begin{enumerate}
\item Using Cayley graphs to construct DSRGs requires us to choose the
correct subset $S$.  Looking at the multiplication table for the set $S$
determines if it works, but what is needed is a way to eliminate
possible subsets before constructing their table.  {\bf Question:} Short of trying every possible subset, is there
a simple way to pick an $S$ that would generate a DSRG, or at least
a better set of criteria so a much smaller set of subsets would need
to be tested?  

\item In using tournaments for all of these constructions, we became
interested in tournaments themselves.  Questions about how many
regular tournaments exist came up, but this is still open and
complete results only exist for small orders.  Also, if a regular
tournament is a circulant matrix, they become easier to work with,
but not all tournaments are circulants.  But, if any regular
tournament is isomorphic to a circulant matrix, which means we can
always choose to work with a circulant representative, our search
can be simplified by allowing us to work matrices with several nice
properties.  A helpful paper on tournaments was \cite{EH}.

\item Since we've shown that isomorphic tournaments, when put through the same construction, will also produce isomorphic graphs, but have found that, for small orders, non-isomorphic tournaments yield non-isomorphic graphs, we would like
to consider how many different classes of non-isomorphic graphs are
created in each case.  {\bf Conjecture: For any construction method $C$, where $C(A)$ is the graph resulting from using matrix $A$ in $C$, $A\cong B$ iff $C(A)\cong C(B)$.}  If this conjecture is true, each parameter set satisfied by our constructions has at least one non-isomorphic graph for each regular tournament, giving us good lower bounds for many parameter sets.  There are a lot of non-isomorphic regular tournaments, especially as the order of the tournaments gets large.  To see exactly how much this gets us, there are already 1223 non-isomorphic regular tournaments of order 11 and this increases greatly to 1495297 non-isomorphic tournaments
of order 13. A complete list of these is found in \cite{Mc}.  If this conjecture is not true, what characteristics of the tournaments are involved in creating isomorphic graphs?  Thinking about how this could happen, the tournaments would have to be similar enough to make isomorphic graphs, yet not be isomorphic themselves, which seems improbable.  However, with as many regular tournaments as there are, it is certainly possible.

\item Almost all of these various constructions are done with regular
tournaments, but this may not be the weakest condition we could
impose to ensure that our constructions still leave us with DSRGs.
There could be a more inclusive class of matrices that also work in
these constructions, allowing for even more non-isomorphic graphs to
be created.

\item The reason that isomorphisms between the resulting graphs is
important is that none of the parameter sets able to be produced by
our constructions are new.  They have all been covered by previous
constructions, most notably J\o rgensen's construction in \cite{Jo}.
So, instead of finding new parameter sets, we've provided several
ways of realizing these parameter sets. 

\item How does isomorphism of tournaments transfer to isomorphism of
graphs through the various constructions?  It doesn't necessarily
carry from one construction to another.  We have one sufficient
condition that guarantees isomorphism, but this may not also be a
necessary condition in this case.  Our last 3 constructions are also
very similar, often making DSRGs with the same parameters.  The
similarity of form in the final adjacency matrix, although not
necessarily in construction of that matrix, may lead to a number of
isomorphic graphs. These block matrix constructions use matrices of the
same form, $$M(B) = \left[ \begin{array}{cc}
B & B^T+I \\
B+I &B^T \end{array} \right], $$ but have different restrictions on
what kind of matrix, $B$, makes up each block.  There is possibly a
more general way to use this form of matrix to construct DSRGs that
includes all of our constructions.  This relates back to the
question of whether having our blocks based on regular tournaments
is the weakest possible condition.  The single form $M(B)$, when
given various types of matrices, creates DSRGs.  Finding more such
forms of block matrices, or finding all possible matrices $B$ so
$M(B)$is a DSRG are possible routes towards discovering more
directed strongly regular graphs.
\end{enumerate}


\section{Lists of Directed Strongly Regular Graphs}
\begin{table}[ht]
\caption{Our Constructions of DSRGs} 
\centering      
\vspace{.1 in}
\begin{tabular}{|p{1 cm}|p{6 cm}|p{4 cm}|p{4 cm}|}  
\hline                        
Lemma & Source & DSRG$(v,k,\mu +1,\mu,\mu)$ & DSRG$(v,k,\mu,\mu-1, \mu)$ \\  [0.5 ex]
\hline
\ref{lem2} & A regular tournament with adjacency matrix $A$ and a DSRG with adjacency matrix $M(A)= \left[ \begin{array}{cc}
A & A^T \\
A & A^T \end{array} \right[ $ or $M(A)= \left[\begin{array}{cc}
A & A \\
A^T & A^T \end{array} \right] $&  & $(4k+2,2k,k,k-1,k)$ \\
\hline \ref{lem8} & A regular tournament with adjacency matrix $A$ and a DSRG with adjacency matrix $M(A)= \left[\begin{array}{cc}
A & A^T +I \\
A+I & A^T \end{array} \right] $ &  & $(4k+2,2k+1,k+1,k,k+1)$ \\
\hline \ref{lem4} & A regular tournament with adjacency matrix $A$ and a DSRG with adjacency matrix $M(A) = \left[ \begin{array}{ccccc}
A & A^T & A & \ldots & A^T \\
A & A^T & A & \ldots & A^T \\
\vdots & \vdots & \vdots & \ddots & \vdots \\
A & A^T & A & \ldots & A^T  \end{array} \right] $ or
$ M(A) = \left[\begin{array}{cccc}
A & A & \ldots & A \\
A^T & A^T & \ldots & A^T \\
\vdots & \vdots  & \ddots & \vdots \\
A^T & A^T & \ldots & A^T   \end{array} \right]$ & & $(w(4k +2), 2wk, wk, w(k-1) wk)$\\
\hline \ref{lem5} & $D(T)$ the adjacency matrix of a doubly regular $(m,2)$ team tournament and a DSRG with adjacency matrix $M(D) = \left[ \begin{array}{cc}
D & D^T +I \\
D+I & D^T \end{array} \right] $ &  (4$m$, 2$m$-1, $m$, $m$-1, $m$-1) where $m\equiv 0\;(mod\;4)$ & \\
\hline \ref{lem6} & A regular tournament with adjacency matrix $A$, DSRG with adjacency matrix $M(D) = \left[ \begin{array}{cc}
D & D^T +I \\
D+I & D^T \end{array} \right] $ &  $(4n+4, 2n+1, n+1, n,n)$ where n is odd &  \\
\hline \ref{lem7} & $L$ is an $2s+2\times 2s+2$ matrix equal to $\Pi +
\Pi^2 + \ldots + \Pi^s$ where $\Pi^n$ is an $n\times n$ matrix as
defined in Lemma \ref{lem3}, DSRG with adjacency matrix $M(L) =
\left[\begin{array}{cc}
L & L^T+I \\
L+I & L^T \end{array} \right] $ &  $(4(s+1), 2s+1, s+1, s, s)$ & \\

\hline
\end{tabular}
\label{DSRGtable}
\end{table}

\begin{table}[ht]
\caption{DSRGs $(v,k,t,\lambda,\mu)$ Constructed by Lemma \ref{lem2}/Lemma \ref{lem4} and DSRGs $(v',k',t',\lambda'\mu')$ Constructed from Lemma \ref{lem8} with up to 34 Vertices} 
\centering      
\vspace{.1 in}
\begin{tabular}{|c|c|c|c|c||c|c|c|c|c|}  
\hline                        
v & k & t & $\lambda$ & $\mu$ & v' & k' & t' & $\lambda$' & $\mu$' \\  [0.5 ex]
\hline
6 & 2 & 1 & 0 & 1& 6 & 3 & 2 & 1 & 2 \\
10 & 4 & 2 & 1 & 2 & 10 & 5 & 3 & 2 & 3\\
14 & 6 & 3 & 2 & 3 & 14 & 7 & 4 & 3 & 4\\
18 & 8 & 4 & 3 & 4 & 18 & 9 & 5 & 4 & 5 \\
22 & 10 & 5 & 4 & 5 & 22 & 11 & 6 & 5 & 6\\
26 & 12 & 6 & 5 & 6 & 26 & 13 & 7 & 6 & 7\\
30 & 14 & 7 & 6 & 7 & 30 & 15 & 8 & 7 & 8\\
34 & 16 & 8 & 7 & 8 & 34 & 17 & 9 & 8 & 9\\

\hline
\end{tabular}
\label{DSRGtable}
\end{table}

\begin{table}[ht]
\caption{DSRGs Constructed from Lemmas \ref{lem5}, \ref{lem6}, and \ref{lem7} with up to 32 Vertices} 
\centering      
\vspace{.1 in}
\begin{tabular}{|p{1 cm}|p{1 cm}|p{1 cm}|p{1 cm}|p{1 cm}||p{4 cm}|p{4 cm}|}  
\hline                        
v & k & t & $\lambda$ & $\mu$  & Consturction Method & Remark on Construction\\  [0.5 ex]
\hline
8 & 3 & 2 & 1 & 1 & Lemmas \ref{lem6} and \ref{lem7} & This graph is unique and therefore Lemmas \ref{lem6} and \ref{lem7} Isomorphic\\
\hline 12 & 5 & 3 & 2 & 2 & Lemma \ref{lem7} & \\
\hline 16 & 7 & 4 & 3 & 3 & Lemmas \ref{lem5}, \ref{lem6}, and \ref{lem7} & Lemmas \ref{lem5} and \ref{lem6} come from one tournament so they are isomorphic, they are not isomorphic to Lemma \ref{lem7} \\
\hline 20 & 9 & 5 & 4 & 4 & Lemma \ref{lem7} \\
\hline 24 & 11 & 6 & 5 & 5 & Lemmas \ref{lem6} and \ref{lem7} & isomorphic (???)\\
\hline 28 & 13 & 7 & 6 & 6 & Lemma \ref{lem7}& \\
\hline 32 & 15 & 8 & 7 & 7 & Lemmas \ref{lem5}, \ref{lem6}, and \ref{lem7} & isomorphic (???)\\

\hline
\end{tabular}
\label{DSRGtable}
\end{table}

\end{document}